\documentclass[12pt]{article}

\usepackage{floatflt}
\usepackage{mathptm}    % use Times fonts for text and math
\usepackage{latexsym}   % use all LaTeX fonts
\usepackage{amssymb}    % use all AMS fonts
\usepackage{amsmath}    % include some AMS-LaTeX functionality
\usepackage{enumerate}

\newtheorem{theorem}{Theorem}[section]
\newtheorem{lemma}[theorem]{Lemma}

\newtheorem{corollary}[theorem]{Corollary}

\newenvironment{prf}[1]{\trivlist
\item[\hskip
\labelsep{\it #1.\hspace*{.3em}}]}{%~\hspace{\fill}~$\square$%
\endtrivlist}
\newenvironment{proof}{\begin{prf}{Proof}}{\end{prf}}

\newtheorem{predefinition}[theorem]{Definition}

\newtheorem{preremark}[theorem]{Remark}
\newenvironment{remark}{\begin{preremark}\rm}{\end{preremark}}
\newtheorem{prenotation}[theorem]{Notation}

\newtheorem{preexample}[theorem]{Example}

\newtheorem{prequestion}[theorem]{Question}

 \makeatletter
\def\emppsubsection{\@startsection{subsection}{2}{\z@}{-3.25ex plus -1ex minus -.2ex}{-1em}}

\makeatother

\newcommand \CH {{\mathcal H}}

\newcommand \PP {{\mathbb P}^1}

\newcommand \ZZ {{\mathbb Z}}

\newcommand  \FF {{\mathbb F}}

\newcommand \GG {{\mathbb G}}

\newcommand \p {{\mathfrak p}}

\newcommand \Aut {\mathop{\rm Aut}}

\newcommand \dime {\mathop{\rm dim}}

\newcommand \Hom {\mathop{\rm Hom}}

\newcommand \Jac {\mathop{\rm Jac}}

\newcommand \nin {\not \in}

\begin{document}

\title{Klein-Four covers of the projective line in characteristic two}
\author{Darren Glass}

\maketitle

\begin{abstract}
In this paper we examine curves defined over a field of
characteristic $2$ which are $(\ZZ/2\ZZ)^2$-covers of the projective
line.  In particular, we prove which $2$-ranks occur for such curves
of a given genus and where possible we give explicit equations for
such curves.
\end{abstract}

\section{Introduction}

There are many ways to stratify the moduli space of curves.  In
characteristic $p>0$, one of the most natural stratifications comes
from looking at the $p$-ranks of the curves. The $p$-rank of a curve
$X$ (or, more precisely, the $p$-rank of its Jacobian) can be
defined as $\dime_{\FF_p} \Hom(\mu_p, \Jac(X))$ where $\mu_p$ is the
kernel of Frobenius on $\GG_m$.  In particular, curves of $p$-rank
$\sigma$ will have precisely $p^\sigma$ distinct $p$-torsion points
on their Jacobian.

There is an idea that curves that have lots of automorphisms should
have small $p$-rank -- in particular, the automorphisms would have
to permute the $p$-torsion points and this would lead to a
restriction on the possible number of these points.  This idea has
never been precisely put into the form of a conjecture or theorem,
but several attempts have been made to investigate the relationship
between automorphism groups and $p$-ranks.

It follows from \cite{GP} in characteristic $p>2$ and \cite{Z} in
characteristic $2$ that hyperelliptic curves behave similar to
generic curves in the sense that there exist curves of each
possible $2$-rank. In this note, we investigate what one can say
about the $2$-ranks of curves which have multiple copies of
$\ZZ/2\ZZ$ in their automorphism group. More precisely, we
consider curves defined over an algebraically closed field of
characteristic $p=2$ which admit an action of $(\ZZ/2\ZZ)^2$ and
such that their quotient by this action is $\PP$.

In Section Two of this paper, we introduce notation and recall some
results from \cite{GP} and \cite{GP2} about the moduli space of
Klein-four covers of the projective line.  We also recall some
results from the theory of Artin-Schreier covers that will be used
to compute the genera and $2$-ranks of the relevant curves. Section
Three is concerned with some nonexistence results, and we prove a
series of results about when various $2$-ranks do not occur.  In the
fourth section, we prove that section three covered all possible
obstructions, and in particular we prove (a stronger version of) the
following theorem.

\begin{theorem}
Let $g \ge 0$ and $0 \le \sigma \le g$. Then there exists a curve
$X$ with $G \cong (\ZZ/2\ZZ)^2 \subseteq \Aut(X)$ and $X/G \cong
\PP$ such that $X$ has genus $g$ and $2$-rank $\sigma$ unless
$\sigma=g-1$ or $g$ is even and $\sigma = 1$.
\end{theorem}

It will follow from the constructions of these curves that they are
all defined over the finite field $\FF_4$ and in most cases they can
be chosen to be defined over $\FF_2$.

We also relate our results to a result of Zhu in \cite{Z}.  In
particular, she proves that there exist hyperelliptic curves of
every possible $2$-rank with no extra automorphisms, while the
following theorem shows precisely when a hyperelliptic curve can
have extra involutions.

\begin{theorem}
There are hyperelliptic curves of genus $g$ and $2$-rank $\sigma$
which contain an additional involution in their automorphism group
if and only if $g \equiv \sigma$ (mod $2$).
\end{theorem}

%Finally, in Section Five we look at $(\ZZ/2\ZZ)^n$-covers of $\PP$
%for $n>2$.  While we do not prove a full classification of which
%$2$-ranks occur for these covers, we prove several results including
%the fact that such covers exist for every genus $g$, unlike the case
%when the characteristic of $k$ is not equal to two.

\section{Notation and Machinery}

In this article, we work over an algebraically closed field $k$ of
characteristic $p=2$.  We wish to examine curves that are
$(\ZZ/2\ZZ)^2$-covers of the projective line $\PP_k$.  In \cite{GP},
we examined such curves defined over algebraically closed fields of
characteristic $p>2$ and in particular we used such curves to
construct hyperelliptic curves with particular group schemes arising
as the $p$-torsion of their Jacobians.  When the characteristic of
$k$ is not equal to two, this Hurwitz space of such covers is
well-defined (for details, see the results of Wewer in
\cite{W:thesis}) and in \cite{GP} we denoted the moduli space of
genus $g$ curves which are $(\ZZ/2\ZZ)^2$-covers of $\PP$ by
$\CH_{g,2}$.   However, when the characteristic of $k$ is equal to
two we are in the situation of wild ramification, and Wewer's
results do not hold.  In particular, it is not clear whether
$\CH_{g,2}$ will be well-defined as a smooth moduli space (see
\cite{Ma} and \cite{Sa} for details on some of the aspects that can
go wrong when defining Hurwitz Spaces associated to wild
ramification). In this note we will abuse notation and define
$\CH_{g,2}$ merely as the set of all $(\ZZ/2\ZZ)^2$-covers of $\PP$
when the characteristic of $k$ equals two.

The key result which we will make use of in this paper is the
following theorem which follows immediately from a result of Kani
and Rosen \cite{KR}.

\begin{theorem}\label{T:KR}
Let $X$ be a curve in $\CH_{g,2}$ and let $H_1$, $H_2$, and
$H_{1,2}$ be the three subgroups of the group $(\ZZ/2\ZZ)^2$ with
respect to a fixed basis.  Furthermore, let $C_1$, $C_2$, and
$C_{1,2}$ be the three quotient curves of $X$ by these subgroups.
Then $$\Jac(X) \sim \prod \Jac(C_S)$$ In particular, if $g_S$ is
the genus of $C_S$ and $\sigma_S$ is the $p$-rank of $C_S$ then we
have that $g_X=g_1+g_2+g_{1,2}$ and
$\sigma_X=\sigma_1+\sigma_2+\sigma_{1,2}$.
\end{theorem}

In the case where the characteristic of $k$ is not equal to $2$,
we were able to show that one could further deduce from the fact
that the degree of this isogeny is a power of two some results
about invariants such as the $a$-number.  While we cannot do this
in the case under consideration in this paper, there is more that
we can say.  In particular, we know that $C_1$ and $C_2$ must be
Artin-Schreier covers, and therefore can be put into the form
$C_i: y^2+y = f_i(x)$ where $f_i$ is a rational function all of
whose poles are of odd order.  In this case, it follows from
results of van der Geer and van der Vlugt in \cite{GV} that the
third quotient is of the form $C_{1,2}:y^2+y=f_1(x)+f_2(x)$.

The following results about the genus and $p$-rank of Artin-Schreier
curves in characteristic two are well known and follow from the
Riemann-Hurwitz and Deuring-Shafarevich formulae, and will be used
throughout this note without reference.

\begin{theorem}
Let $y^2+y=f(x)$ define a hyperelliptic curve $C$ in
characteristic two.  Let $f(x)$ have $k$ poles given by
$x_1,\ldots,x_k$ and let $n_i$ be the order of the pole at $x_i$.
Without loss of generality we can assume that all of the $n_i$ are
odd.  Then the genus of $C$ is given by the formula
$-1+\frac{1}{2}\sum(n_i+1)$ and the $2$-rank of $C$ is given by
$k-1$.
\end{theorem}

To conclude this introduction we define the type of a curve $X \in
\CH_{g,2}$ to be the unordered triple $\{g_1,g_2,g_3\}$ consisting
of the genera of the three $\ZZ/2\ZZ$ quotients of $X$.  In
particular, it follows that the $g_i$ are integers such that $0 \le
g_i \le \frac{g+1}{2}$ and $g_1+g_2+g_3=g$. In \cite{GP2} we show
that, when the characteristic of $k$ is {\it not} equal to two, the
irreducible components of $\CH_{g,2}$ correspond to the set of
curves of a given type.  However, as discussed above, when the
characteristic of $k$ is equal to two, the objects $\CH_{g,2}$ may
not be well-defined as geometric objects.  For a given partition
$\p$ of $g$ satisfying the necessary conditions we will again abuse
notation and define $\CH_{g,2,\p}$ to be the set of all Klein-four
covers of $\PP$ of type $\p$.

We note that the type of a curve is technically the type of the
cover $X \rightarrow \PP$, and in a small number of cases a curve
$X$ can be considered a $(\ZZ/2\ZZ)^2$-cover of $\PP$ in more than
one way leading to different types.  We show in \cite{GP2} that this
is rare in characteristic $p \ne 2$ (and happens exactly in the case
where $1 \in \p$).  Similar results hold if $p=2$.

\section{Nonexistence Results}

\begin{theorem}
There are no almost-ordinary curves (curves with $2$-rank
$\sigma=g-1$) in $\CH_{g,2}$ for any $g$.
\end{theorem}

\begin{proof}
Let $X$ be a curve in $\CH_{g,2}$ which is almost-ordinary. It
follows from Theorem \ref{T:KR} that one of its $\ZZ/2\ZZ$
quotients must be almost-ordinary and the other two must be
ordinary.  Let $C_1$ and $C_2$ be the two quotients which are
ordinary so that $C_1$ (resp. $C_2$) is defined by the equation
$y^2+y=f_1(x)$ (resp. $f_2(x)$) where $f_1$ (resp. $f_2$) only has
simple poles. Then $f_1+f_2$ must also have only simple poles and
therefore the curve $C_{1,2}$, which is defined by
$y^2+y=f_1(x)+f_2(x)$, must also be ordinary. This gives a
contradiction.
\end{proof}

In some cases it happens that a given $2$-rank can occur for curves
of some type in $\CH_{g,2}$ but not for curves of another type, as
the following theorems indicate.

\begin{theorem}
There exist curves in $\CH_{g,2,\p}$ of $2$-rank zero only if the
two largest entries of $\p$ are identical. (ie. the elements of
$\p$ can be put into order $g_1=g_2 \ge g_3$)
\end{theorem}

\begin{proof}
Let $X \in \CH_{g,2}$ be a curve whose $2$-rank equals zero.  Then
we can conclude from Theorem \ref{T:KR} that all three of the
hyperelliptic quotients have $2$-rank zero and therefore they can
each be defined by $y^2+y=f_i(x)$ where each $f_i$ has a single
pole.  It follows that (at least) two of these three functions
must have a pole of the same order, and therefore that (at least)
two of the associated subgenera must be identical.
\end{proof}

We define a partition $\p$ to be unbalanced if it contains an
element which is at least $\frac{g}{2}$.  In particular, unbalanced
partitions are of the form $\{\frac{g}{2},*,*\}$ or
$\{\frac{g+1}{2},*,*\}$ depending on the parity of $g$.

\begin{theorem}\label{T:prank1}
There exist smooth curves in $\CH_{g,2,\p}$ with $2$-rank equal to
one if and only if $g$ is odd and $\p$ is unbalanced. In
particular, if $g$ is even then there are no curves in $\CH_{g,2}$
with $2$-rank equal to one.
\end{theorem}

\begin{proof}
Let $X \in \CH_{g,2}$ be a curve whose $2$-rank equals $1$.  Then
 two of the hyperelliptic quotients must have
$2$-rank zero while the third has $2$-rank one.  It follows
without loss of generality that $f_1$ has a pole of order $a$ at
one point and $f_2$ has a pole of order $b$ at another point where
$a$ and $b$ are both odd. In that case we can compute that the
curve $X$ is of type
$\{\frac{a-1}{2},\frac{b-1}{2},\frac{a+b}{2}\}$ which in turn
implies that the genus of the curve $X$ is $a+b-1$ (and is thus
odd) while the genus of the curve $C_3$ is $\frac{a+b}{2} =
\frac{g+1}{2}$.
\end{proof}

A quite different but similarly restrictive results holds if we
look at curves with $2$-rank equal to $2$.

\begin{theorem}\label{T:prank2}
Let $\p = \{g_1,g_1,g_1\}$ be a totally balanced partition.  Then
there do not exist curves in $\CH_{g,2,\p}$ with $2$-rank equal to
two.
\end{theorem}

\begin{proof}
Let $X \in \CH_{g,2}$ be a curve whose $2$-rank is equal to $2$.
Let $C_1,C_2$, and $C_3$ be the three quotient curves and let
$\sigma_i$ be the $2$-rank of $C_i$.  Then it follows from Theorem
\ref{T:KR} that (without loss of generality) either $\sigma_1 = 2$
and $\sigma_2=\sigma_3=0$ or $\sigma_1=\sigma_2=1$ and
$\sigma_3=0$. However, the first case cannot happen, because it
would imply that $f_1$ would have $3$ poles while each of $f_2$
and $f_3$ would have a single pole.

Therefore we must be in the second case, in which $f_1$ and $f_2$
each have two poles and $f_3$ has one pole. It follows that we can
assume that $f_1$ and $f_2$ each have poles at zero and infinity
and that $f_3$ has a pole only at infinity. Without loss of
generality, we may assume that $ord_\infty(f_1) \ge
ord_\infty(f_3)$ which will in turn imply that $g_1>g_3$.
Therefore, $\p$ cannot be a totally balanced partition.
\end{proof}

\begin{theorem} \label{T:unbalanced}
Let $\p$ be an unbalanced partition.  Then there exist smooth
curves in $\CH_{g,2,\p}$ with $2$-rank equal to $\sigma$ only if
$g \equiv \sigma$ (mod $2$).
\end{theorem}

\begin{proof}
Assume that $X$ is a curve in $\CH_{g,2,\p}$ and that
$\frac{g+1}{2} \in \p$, so that $g$ must be odd.  Then there
exists an involution $\tau \in Aut(X)$ such that the genus of
$C_1=X/<\tau>$ is equal to $\frac{g+1}{2}$.  It follows from the
Riemann-Hurwitz formula that this cover must be etale.  Therefore,
if we apply the Deuring-Shafarevich formula we see that $\sigma_X
= 2 \sigma_{C_1}-1$ is odd.

Similarly, if we assume that $X$ is a curve in $\CH_{g,2,\p}$ and
that $\frac{g}{2} \in \p$ (and hence $g$ is even), we note that it
follows from the Riemann-Hurwitz formula that the cover $X \to
C_1$ must be ramified at a single point.  Again, it will follow
from Deuring-Shafarevich that $\sigma_X = 2 \sigma_{C_1}$ must be
even.

Therefore, in both cases where we look at curves whose associated
partitions are unbalanced we see that $\sigma_X \equiv g_X$ (mod
$2$). The converse will follow from results in Section
\ref{S:exist}
\end{proof}

\section{Existence Results} \label{S:exist}

The main result in this section is that all $2$-ranks occur for
all types of curves in $\CH_{g,2}$ except for those that we have
proven do not occur in the previous section.  In particular, we
will prove the following theorem.

\begin{theorem}\label{T:main}
There exist curves of $2$-rank $\sigma$ in $\CH_{g,2,\p}$ for all
$0 \le \sigma \le g$ except in the following cases:

\begin{enumerate}[i.]
\item $\sigma = 0$, $\p \ne \{g_1,g_1,g_3\}$ with $g_3 \le g_1$.
\item $\sigma = 1$, $\frac{g+1}{2} \nin \p$.
\item $\sigma = 2$, $\p$ totally balanced.
\item $\sigma = g-1$.
\item $\sigma \not \equiv g$ (mod $2$), $\p$ unbalanced.
\end{enumerate}

\end{theorem}

We will prove this theorem by induction on $\sigma$ after looking
at some base cases.  Throughout these, let $\alpha$ be one of the
elements of $\FF_4$ other than one or zero.

\begin{lemma}\label{L2}
Let $\p$ be a partition which is neither completely balanced or,
if $g$ is odd, unbalanced. Then there exist curves in
$\CH_{g,2,\p}$ of $2$-rank equal to two.
\end{lemma}

\begin{proof}
Let $\p = \{g_1,g_2,g_3\}$ with $g_1 \ge g_2 \ge g_3$.  Let $a =
2g_3+1$, $b=2(g_1-g_3)-1$ and $c=2(g_2+g_3-g_1)+1$.  It is clear
that $a,b,$ and $c$ are all odd, and that $a \ge c$.  Furthermore,
$b \ge 1$ because $\p$ is not completely balanced and $c\ge 1$
because $g_1 \le g/2$.  If we now let $f_1 = x^a + \frac{1}{x^b}$
and $f_2 = \alpha x^c+\frac{1}{x^b}$ we see that
$f_3=f_1+f_2=x^a+\alpha x^c$ and a simple computation shows that
the fibre product $X$ will have $2$-rank equal to $2$ and lie in
$\CH_{g,2,\p}$.
\end{proof}

\begin{lemma}\label{L3}
If $\sigma = 3$ then Theorem \ref{T:main} holds.  In particular,
there are curves of $2$-rank equal to three in each $\CH_{g,2,\p}$
except the case where $g$ is even and $\frac{g}{2} \in \p$.
\end{lemma}

\begin{proof}
We consider two different cases.  First, if $g$ is odd and $\p$ is
unbalanced then $\frac{g+1}{2} \in \p$ and in particular we may
assume that the branch loci of $f_1$ and $f_2$ will be disjoint.
If we let $f_1 = x^a + \frac{1}{x-1}$ and $f_2 = \frac{1}{x^b}$
then we will have $\sigma_1=1, \sigma_2=0,$ and $\sigma_3 = 2$ so
that $\sigma_X = 3$. Furthermore, $g_1=\frac{a+1}{2}$,
$g_2=\frac{b-1}{2}$, and $g_3=\frac{a+b+2}{2} = g_1+g_2+1 =
\frac{g+1}{2}$.  We can choose $a$ and $b$ in order to get any
unbalanced partition of $g$ that we desire.

For the other case, assume that $\p = \{g_1,g_2,g_3\}$ with $1 \le
g_3 \le g_2 \le g_1 \le \frac{g-1}{2}$.  Set $a=2(g_1-g_3)+1$,
$b=2g_3-1$ and $d=2(g_2+g_3-g_1)-1$.  We note that our hypotheses
imply that $a$, $b$, and $d$ are all odd positive numbers with $b
\ge d$.  Now, let $f_1 = x^a+ \frac{\alpha}{x^b}$, $f_2 = x^a +
\frac{1}{x^d}$ and $f_3=f_1+f_2$. Then the curve defined by
$y^2+y=f_i$ will have genus $g_i$ and $2$-rank equal to one, and
therefore $X$ will be a curve of $2$-rank equal to three in
$\CH_{g,2,\p}$.
\end{proof}

\begin{lemma}\label{L4}
If $\sigma = 4$ then Theorem \ref{T:main} holds.  In particular,
there are curves of $2$-rank equal to four in each $\CH_{g,2,\p}$
except the case where $g$ is odd and $\frac{g+1}{2} \in \p$.
\end{lemma}

\begin{proof}
Assume that $\p=\{g_1,g_2,g_3\}$ where $g_1 > g_2 \ge g_3$. Let
$a=2g_2-1$, $b=2(g_1-g_2)-1$, and $c=2(g_2+g_3-g_1)+1$.  One can
easily check that $a,b,$ and $c$ are all positive odd numbers
(recall that the fact that $\frac{g+1}{2} \nin \p$ implies that $g_1
\le g_2+g_3$) and furthermore that $a \ge c$.  Let $f_1 = x^a +
\frac{1}{x^b} + \frac{1}{x+1}$, $f_3=x^c+\frac{1}{x^b}$, and
$f_2=f_1+f_3$. Then the curve defined by the equation $y^2+y=f_i(x)$
has genus $g_i$ and the fibre product will have $2$-rank $\sigma=4$
as desired.

On the other hand, assume that $g_1=g_2 \ge g_3 \ge 2$.  In this
case, let $a=2g_1-1$ and $b=2g_3-3$.  Then it is clear that $a$ and
$b$ are positive odd integers with $a > b$.  If we define
$f_1=x^a+\frac{1}{x}$ and $f_3=x^b+\frac{1}{x}+\frac{1}{x+1}$ we can
see that the curves will have the desired properties.

It remains to consider the partitions
$\p=\{\frac{g}{2},\frac{g}{2},0\}$ in the case where $g$ is even and
$\p=\{\frac{g-1}{2},\frac{g-1}{2},1\}$ if $g$ is odd.  For the
former, we let $f_1=x^{g-3}+\frac{1}{x}+\frac{1}{x+1}$ and $f_2 =
\alpha x$ and $f_3=f_1+f_2$.  Then $g_1=g_3=g/2$ and $g_2=0$ while
$\sigma_1=\sigma_3=2$ and $\sigma_2=0$.  For the latter, we note
that $g$ must be odd and $g>\sigma+1=5$, so we may assume that $g
\ge 7$.  Let $f_1=x^{g-4}+\frac{1}{x}+\frac{1}{x+1}$ and $f_2 =
\alpha x$.  These equations define curves with the desired genera
and $2$-ranks.

\end{proof}

\begin{remark}\label{R:induct}
Now that we have shown that the theorem is true for $\sigma \le 4$
we are ready to consider the inductive step.  The key idea is to
notice that if there is a curve $X$ in $\CH_{g,2,\{g_1,g_2,g_3\}}$
with $2$-rank equal to $\sigma$ then there will be a curve
$\tilde{X}$ in $\CH_{g+3,2,\{g_1+1,g_2+1,g_3+1\}}$ with $2$-rank
equal to $\sigma+3$.  In particular, if the three hyperelliptic
quotients of $X$ are defined by the equations $y^2+y=f_i(x)$, then
without loss of generality we may assume that none of the $f_i$ have
poles at infinity.  Then we define $\tilde{f_1}=f_1+ x$,
$\tilde{f_2} = f_2 + \alpha x$ and $\tilde{f_3} = f_3 +
(\alpha+1)x$.  It is clear that
$\tilde{f_3}=\tilde{f_1}+\tilde{f_2}$ and that the curve $\tilde{X}$
defined by the fibre product of $y^2+y=\tilde{f_1}(x)$ and
$y^2+y=\tilde{f_2}(x)$ will lie in
$\CH_{g+3,2,\{g_1+1,g_2+1,g_3+1\}}$ and have $2$-rank $\sigma +3$.
Therefore, once we have shown that Theorem \ref{T:main} holds for
$\sigma$ we have {\it almost} shown that it will hold for
$\sigma+3$.  However, to be complete there are still a few cases we
must consider.
\end{remark}

\begin{lemma}\label{L:unbeven}
The component $\CH_{g,2,\{g/2,g/2,0\}}$ contains curves of every
even $2$-rank.
\end{lemma}

\begin{proof}
In order to construct a curve in $\CH_{g,2,\{g/2,g/2,0\}}$ with
$2$-rank equal to $2k$ we first note that we can find a
hyperelliptic curve of genus $g/2$ with $2$-rank equal to $k$. Let
us assume that this curve is defined by the equation
$y^2+y=f_1(x)$ where $f_1$ has a pole at infinity.  Let $f_2$ be
some constant multiple of $x$ so that $f_3=f_1+f_2$ will have the
same poles (with the same orders) as $f_1$.  Then it follows from
our construction that the curve $X$ will have $2$-rank $2k$ and
will lie in $\CH_{g,2,\{g/2,g/2,0\}}$.
\end{proof}

\begin{lemma}\label{L:unbodd}
Let $g$ be odd and $\p$ be an unbalanced partition of $g$.  Then
$\CH_{g,2,\p}$ contains curves of all odd $p$-ranks $\sigma =
2k+1$.
\end{lemma}

\begin{proof}
If $g$ is odd and $\p$ is unbalanced, then $\p =
\{\frac{g+1}{2},g_1,g_2\}$ where $g_1 \ge g_2$.  We note that we
can construct hyperelliptic curves $C_1$ and $C_2$ so that the
genus of $C_i$ is $g_i$ and the $2$-rank of $C_i$ is $k_i$ for all
$0 \le k_i \le g_i$.  Furthermore, we can assume that the branch
loci are distinct.  If we let $X$ be the fibre product of $C_1$
and $C_2$ and consider the third hyperelliptic quotient of $X$ we
see that it will have genus $g_1+g_2+1$ and $2$-rank $k_1+k_2+1$.
If we choose $k_1$ and $k_2$ so that $k_1+k_2=k$ then $X$ will
have $2$-rank equal to $\sigma$ and lie in $\CH_{g,2,\p}$.
\end{proof}

\begin{lemma}\label{L:g-1/2}
Let $g$ be odd and $\frac{g-1}{2} \in \p$ but $\frac{g+1}{2} \nin
\p$. Then there are curves in $\CH_{g,2,\p}$ with $2$-rank equal to
$2k$ for all $0 \le k \le \frac{g-3}{2}$.
\end{lemma}

\begin{proof}
Let $\p = \{\frac{g-1}{2},g_1,g_2\}$ with $g_1 \ge g_2>0$ and let
$\sigma=2k$ be as above. Because $\sigma \le g-3$ we have that $k
\le g_1+g_2-2$ and therefore we can choose $k_1$ and $k_2$ so that
$k_1+k_2=k$ but $k_i < g_i$. In particular, we can define a
function $h_1(x)$ which has $k_1$ poles (none of which are at
infinity) so that if we look at $\frac{1}{2}\sum(n_i+1)$ where the
sum runs over the poles of $h_1$, each of which is of order $n_i$,
then we can get any integer which is at least $k_1$ and in
particular we can get $g_1-1$.  To be precise, we can choose $h_1$
so that the curve $C_1$ defined by $y^2+y=x^3+h_1(x)$ will have
genus $g_1$ and $2$-rank $k_1$.  Similarly, we can choose $h_2$
with poles distinct from those of $h_1$ so that the curve $C_2$
defined by $y^2+y=\alpha x^3 + h_2(x)$ will have genus $g_2$ and
$2$-rank $k_2$.

If we look at the normalization of the fibre product of $C_1$ and
$C_2$ we see that the third quotient will be defined by the
equation $y^2+y = (\alpha + 1) x^3 + h_1(x) +h_2 (x)$ and
therefore will have genus $g_1+g_2-1$ and $2$-rank $k_1+k_2=k$.
Thus, the curve $X$ lies in $\CH_{g,2,\{\frac{g-1}{2},g_1,g_2\}}$
and has $2$-rank equal to $2k$, as desired.
\end{proof}

Before we prove the main theorem in general, we look at the case
where $\sigma=5$, in which we need to fill an extra gap.

\begin{lemma}
If $\sigma = 5$ then Theorem \ref{T:main} holds.  In particular, if
$g \ge 7$ there are curves of $2$-rank equal to five in each
$\CH_{g,2,\p}$ except the case where $g$ is even and $\frac{g}{2}
\in \p$.
\end{lemma}

\begin{proof}
Let $\p = \{g_1,g_2,g_3\}$ be a partition of $g$ with $0 \le g_3
\le g_2 \le g_1 \le \frac{g+1}{2}$.  We wish to show that there
are curves in $\CH_{g,2,\p}$ of $2$-rank equal to five unless
$g_1=\frac{g}{2}$ (in which case $g$ will be even).  If $g_1 =
\frac{g+1}{2}$ then the result follows from Lemma \ref{L:unbodd}.

If $g_1 \le \frac{g-1}{2}$ then it follows that $g_3>0$ and thus
$\hat{\p} = \{g_1-1,g_2-1,g_3-1\}$ gives a partition of $g-3$ all
of whose entries are at most $\frac{g-3}{2}$.  Thus, by Lemma
\ref{L2} there are curves in $\CH_{g-3,2,\hat{p}}$ of $2$-rank
equal to $2$ unless $\hat{p}$ (and therefore $\p$ is completely
balanced).  By the induction argument in Remark \ref{R:induct} we
therefore have curves whose $2$-rank is equal to five in
$\CH_{g,2,\p}$.

It remains to consider the case where $\p$ is totally balanced:
that is, where $g_1=g_2=g_3=a$.  To deal with this case, let
$f_1=x^a + \frac{1}{x^a}$ and $f_2 = x^a + \frac{1}{(x-1)^{a-2}} +
\frac{1}{x-\alpha}$ and $f_3=f_1+f_2$.  One can easily compute
that these choices will lead to a curve $X \ in
\CH_{g,2,\{a,a,a\}}$ whose $2$-rank is equal to $5$.
\end{proof}

We are finally ready to prove Theorem \ref{T:main}.

\begin{proof}
Given the results of the above lemmata, it suffices to consider
the case where $\sigma \ge 6$. In this case, we only need to prove
that there are curves of $2$-rank equal to $\sigma$ in every
partition if $g \equiv \sigma$ (mod $2$) and that there are curves
of $2$ rank equal to $\sigma$ in every partition whose entries are
all at most $\frac{g-1}{2}$ if $g \not \equiv \sigma$ (mod $2$).

If $0 \in \p$ then $\p$ must be unbalanced, and therefore we only
need to consider the case where $g \equiv \sigma$ (mod $2$).  The
result then follows from Lemma \ref{L:unbodd} if $g$ is odd and
from Lemma \ref{L:unbeven} if $g$ is even.  Similarly, if
$\frac{g+1}{2} \in \p$ the result follows from Lemma
\ref{L:unbodd}.

Next we consider the case where $\p=\{\frac{g}{2},g_1,g_2\}$ with
$g_1$ and $g_2$ both positive.  It follows from Theorem
\ref{T:unbalanced} that there are no curves in $\CH_{g,2,\p}$
whose $2$-rank is odd, so we wish to show that there will be
curves of all even $2$-ranks.  We note that $\hat{p} =
\{\frac{g}{2}-1,g_1-1,g_2-1\}$ gives a partition of $g-3$.  If $g
\equiv \sigma$ (mod $2$) then $g-3 \equiv \sigma-3$ and there will
be curves in $\CH_{g-3,2,\hat{\p}}$ of $2$-rank $\sigma-3$ by the
inductive hypotheses which can be used to construct curves in
$\CH_{g,2,\p}$ of $2$-rank equal to $\sigma$ by the inductive
argument in \ref{R:induct}.

If $\p = \{\frac{g-1}{2},g_1,g_2\}$ we note that both $g_1$ and
$g_2$ must be positive. It follows from Lemma \ref{L:g-1/2} that
there are curves of every even $2$-rank less than $g-1$ in
$\CH_{g,2,\p}$. To construct the curves of odd $2$-rank $\sigma$,
we note that $\hat{\p} = \{\frac{g-3}{2},g_1-1,g_2-1\}$ gives a
partition of $g-3$, and that $g-3 \equiv \sigma-3$ (mod $2$) and
therefore there are curves in $\CH_{g-3,2,\hat{p}}$ of $2$-rank
equal to $\sigma-3$.  The result then follows from the inductive
process described in Remark \ref{R:induct}.

If all entries of $\p$ are at least $1$ and at most
$\frac{g-2}{2}$, we note $\hat{\p} = \{g_1-1,g_2-1,g_3-1\}$ gives
a partition of $\hat{g} = g-3$ such that each $\hat{g_i} = g_i-1$
is at most $\frac{\hat{g}-1}{2}$ and therefore there exist curves
of $2$-rank $\sigma-3$ in $\CH_{\hat{g},2,\hat{\p}}$.  By the
inductive procedure described in Remark \ref{R:induct} we can
construct a curve in $\CH_{g,2,\p}$ with $2$-rank equal to
$\sigma$, proving the theorem.

\end{proof}

In \cite{Z}, Zhu proves that there exist hyperelliptic curves with
no extra automorphisms of every possible $2$-rank.  The following
result shows that there are often hyperelliptic curves that {\it
do} have an extra involution.

\begin{corollary}
There are hyperelliptic curves of genus $g$ and $2$-rank $\sigma$
which contain an additional involution in their automorphism group
if and only if $g \equiv \sigma$ (mod $2$).
\end{corollary}

\begin{proof}
It is well known that the hyperelliptic involution is contained in
the center of the automorphism group of a curve.  Therefore, if
there is another involution in the automorphism group then we must
have a Klein-four action on the curve and therefore we will be in
the setup above.  Furthermore, it follows that the partition $\p$
corresponding to this curve contains a zero and is therefore
either $\p = \{\frac{g+1}{2},\frac{g-1}{2},0\}$ or
$\p=\{g/2,g/2,0\}$.  In either case, the partition is unbalanced
and therefore $g \equiv \sigma$ (mod $2$) by Theorem
\ref{T:unbalanced}.

Conversely, it follows from Theorem \ref{T:main} that if $g \equiv
\sigma$ (mod $2$) then there will exist curves in this partition,
which will therefore be both hyperelliptic and contain an extra
involution.
\end{proof}

We note that this does not answer the question of the automorphism
groups fully, as the curves may have automorphisms of degree greater
than two.  We examine the question of the possible $2$-ranks of
hyperelliptic curves with extra automorphisms in depth in \cite{G2}.

\bibliographystyle{abbrv}
\bibliography{prank}
\end{document}